\documentclass [12pt] {amsart}
\usepackage{times}
\usepackage{latexsym,amssymb, amsmath}
\usepackage{enumerate}
\usepackage[margin=2.75cm]{geometry}

\def\pf{\noindent\emph{Proof: }}

\usepackage{amsmath, amsfonts, amssymb, amsthm}
\newtheorem{thm}{Theorem}
\newtheorem{cor}[thm]{Corollary}
\newtheorem{lemma}[thm]{Lemma}

\numberwithin{thm}{section}

\begin{document}

\title{A Hopf's Lemma for Higher Order Differential Inequalities and Its Applications}

\author{Yifei Pan}
\address{Department of Mathematical Sciences\\
 Indiana University-Purdue University Fort Wayne\\
 Fort Wayne, Indiana 46805} 
\email{pan@ipfw.edu}

\author {Mei Wang}
\address {Department of Statistics\\University of Chicago \\Chicago, Illinois 60637}
\email{meiwang@galton.uchicago.edu}

\author { Yu Yan }
\address {Department of Mathematics and Computer Science\\ Huntington University\\Huntington, Indiana 46750}
\email {yyan@huntington.edu}

\begin{abstract}
We establish a sequential Hopf's Lemma for higher order differential inequalities in one variable and give some applications of this result.  
\end{abstract}

\maketitle

\newtheorem{Thm}{Theorem}
\newtheorem{Lemm}{Lemma}
\newtheorem{Cor}{Corollary}

\section {Introduction}

The Hopf's Lemma is one of the fundamental tools in the study of elliptic partial differential equations \cite{GT}.  There have been many variations and generalizations of this lemma, for example \cite{Li-Ni_2}, \cite{Li-Ni_1}, and \cite{Li-Ni_3}.  But there appears to be no work in the literature on the Hopf's lemma for third or higher order equations, perhaps partially because the maximum principle fails for higher order equations.

\vspace{.1in}
In this paper we study this question in the one dimensional case and prove a sequential Hopf's lemma of higher order in one variable.  One application of this result is the following comparison theorem for $n$th order nonlinear differential operators.

\begin{thm}
\label{thm:functional}
Assume that $K(z_1,... , z_{n+2}   )$ is Lipschitz in all variables and $\frac{\partial K}{\partial z_{n+2} }>0$ almost everywhere, where $n \geq 2$.   Suppose $u(x)$ and $v(x)$ are two functions in $C^{n}((a,b))$ that satisfy
\begin{equation}
\label{eq:functional}
K \left ( x, u(x), u'(x), ... , u^{(n)}(x) \right )  \leq K \left ( x, v(x), v'(x), ... , v^{(n)}(x) \right )   \quad \quad \text{for all $x \in (a,b)$}
\end{equation} 
and 
\begin{equation*}
u(x_0)=v(x_0), \quad u'(x_0)=v'(x_0), \quad ...\,\, , \quad u^{(n-1)}(x_0)=v^{(n-1)}(x_0)  \quad \quad  \text{for some $x_0 \in (a,b)$} .
\end{equation*}

\vspace{.05in}
\noindent
If $n$ is even, then there exists $\delta >0$ such that $u(x) \leq v(x)$ for $x \in (x_0 - \delta, x_0 + \delta) .$

\vspace{.05in}
\noindent
If $n$ is odd, then there exists $\delta >0$ such that $u(x) \geq v(x)$ for $x \in (x_0 - \delta, x_0) $ and $u(x) \leq v(x)$ for $x \in (x_0 , x_0 + \delta) $.

\end{thm}

\vspace{.05in}
This theorem shows that if $u$ and $v$ have $(n-1)$-th order of contact at a point $x_0$, then they intersect only once in a small neighborhood of $x_0$.  The crucial ingredient in the proof is a higher order sequential version of Hopf's lemma. 

\begin{thm}
\label{thm:Hopf}
Let $u \in C^{n}((a,b))\bigcap C^{n-1}([a,b))$ be a function which satisfies  
\begin{equation}
\label{eq:main}
u^{(n)}(x)+a_{n-1}(x)u^{n-1}(x)+\cdots +a_1(x)u'(x)+a_0(x)u(x) \leq 0  \hspace{.2in} \text{for } \,\, x\in (a,b),
\end{equation}
 where $n \geq 2$ is a positive integer and $a_{n-1}(x), ... , a_1(x), a_0(x) $ are in $C([a,b))$. Suppose $u$ satisfies 
\begin{equation}
\label{eq:initial}
u(a) \,\, = \,\, u'(a) \,\, = \,\, \cdots  \,\, = u^{(n-2)}(a) \,\, = \,\, 0 ,
\end{equation}
and 
\begin{equation}
\label{eq:sequence}
\text{there exists a sequence } \{x_i\} \,\, \text {such that} \,\, a<x_i<b ,\,\, x_i \to a, \,\, \text{ and } \,\, u(x_i)>0 .
\end{equation}

\noindent
Then $u^{(n-1)}(a)>0$.  Furthermore, $u>0$ in a neighborhood of $a$.

\vspace{.1in}
\noindent
When $n=2$, it suffices to assume that $a_1(x)$ and $a_0(x) $ are bounded functions.

\end{thm}

The Taylor's expansion of $u$ at $a$ and Condition (\ref{eq:sequence}) easily imply that $u^{(n-1)}(a) \geq 0$, so the key is that it is strictly positive.  At the right side endpoint of an interval, we have 

\begin{thm}
\label{thm:Hopf_right}
Let $u \in C^{n}((a,b))\bigcap C^{n-1}((a,b])$ be a function which satisfies  
\begin{equation*}
\label{eq:main_right}
u^{(n)}(x)+a_{n-1}(x)u^{n-1}(x)+\cdots +a_1(x)u'(x)+a_0(x)u(x) \leq 0  \hspace{.2in} \text{for } \,\, x\in (a,b),
\end{equation*}
 where $n \geq 2$ is a positive integer and $a_{n-1}(x), ... , a_1(x), a_0(x) $ are in $C((a,b])$. Suppose $u$ satisfies 
\begin{equation}
\label{eq:initial_right}
u(b) \,\, = \,\, u'(b) \,\, = \,\, \cdots  \,\, = u^{(n-2)}(b) \,\, = \,\, 0 .
\end{equation}

\vspace{.05in}
 
\noindent
If $n$ is even and
\begin{equation}
\label{eq:sequence_even}
\text{there exists a sequence } \{x_i\} \,\, \text {such that} \,\, a<x_i<b ,\,\, x_i \to b, \,\, \text{ and } \,\, u(x_i)>0 ,
\end{equation}
\vspace{.05in}
\noindent
then $u^{(n-1)}(b)<0$ and $u>0$ in a neighborhood of  $b$.

\vspace{.05in}

\noindent
If $n$ is odd and
\begin{equation}
\label{eq:sequence_odd}
\text{there exists a sequence } \{x_i\} \,\, \text {such that} \,\, a<x_i<b ,\,\, x_i \to b, \,\, \text{ and } \,\, u(x_i)<0 ,
\end{equation}
\vspace{.05in}
\noindent
then $u^{(n-1)}(b)<0$ and $u<0$ in a neighborhood of  $b$.

\vspace{.1in}
\noindent
When $n=2$, it suffices to assume that $a_1(x)$ and $a_0(x) $ are bounded functions.

\end{thm}

\vspace{.1in}
In the subsequent sections we will prove the above theorems and discuss some applications.

\vspace{.2in}

\section{Proof of the Comparison Theorem}
\label{sec:comparison}

\vspace{.05in}
\noindent
Since $u$ must be negative or $0$ near $a$ if condition (\ref{eq:sequence}) is not met, an equivalent statement of Theorem \ref{thm:Hopf} is

\vspace{.05in}

\begin{thm}
\label{thm:Hopf-equivalent}
Let $u \in C^{n}((a,b)) \bigcap C^{n-1}([a,b)) $ be a function which satisfies  
\begin{equation*}
\label{eq:main_equivalent}
u^{(n)}(x)+a_{n-1}(x)u^{n-1}(x)+\cdots +a_1(x)u'(x)+a_0(x)u(x) \leq 0 \hspace{.2in} \text{for } \,\, x\in (a,b),
\end{equation*}
 where $n\geq 2$ is a positive integer and $a_{n-1}(x), ... , a_1(x), a_0(x) $ are in $C([a,b))$. 

\noindent 
If 
\begin{equation*}
\label{eq:initial_equivalent}
u(a) \,\, = \,\, u'(a) \,\, = \,\, \cdots  \,\, = u^{(n-2)}(a) \,\, = u^{(n-1)}(a) \,\, =\,\, 0 ,
\end{equation*}

\vspace{.05in}
\noindent
then $u(x) \leq 0$ for all $x$ sufficiently close to $a$.

\vspace{.05in}
\noindent
When $n=2$, it suffices to assume that $a_1(x)$ and $a_0(x) $ are bounded functions.

\end{thm}

\noindent
Similarly, an equivalence of Theorem \ref{thm:Hopf_right} is

\vspace{.05in}

\begin{thm}
\label{thm:Hopf-equivalent_right}
Let $u \in C^{n}((a,b)) \bigcap C^{n-1}((a,b]) $ be a function which satisfies  
\begin{equation*}
\label{eq:main_equivalent}
u^{(n)}(x)+a_{n-1}(x)u^{n-1}(x)+\cdots +a_1(x)u'(x)+a_0(x)u(x) \leq 0 \hspace{.2in} \text{for } \,\, x\in (a,b),
\end{equation*}
 where $n\geq 2$ is a positive integer and $a_{n-1}(x), ... , a_1(x), a_0(x) $ are in $C((a,b])$. 
 
\noindent 
Suppose
 $$
 u(b) \,\, = \,\, u'(b) \,\, = \,\, \cdots  \,\, = u^{(n-2)}(b) \,\, = u^{(n-1)}(b) \,\, =\,\, 0 .
$$

\vspace{.1in} 
\noindent
If $n$ is even, then $u(x) \leq 0$ for all $x$ sufficiently close to $b$.

\vspace{.1in} 
\noindent 
If $n$ is odd, then $u(x) \geq 0$ for all $x$ sufficiently close to $b$.

\vspace{.1in}
\noindent
When $n=2$, it suffices to assume that $a_1(x)$ and $a_0(x) $ are bounded functions.

\end{thm}

\vspace{.05in}

\noindent   
Note that
\begin{eqnarray*}
&   &  K \left ( x, u(x), u'(x), ... , u^{(n)}(x) \right )  - K \left ( x, v(x), v'(x), ... , v^{(n)}(x) \right )   \\
& = &  c_0(x)(u-v)+ \cdots + c_{n-1}(u^{(n-1)}-v^{(n-1)} )+ c_{n}(u^{(n)}-v^{(n)})
\end{eqnarray*}

\noindent
where 
\begin{eqnarray*}
c_0(x) & = & \int _0 ^1 \frac{\partial K}{\partial z_2} \left (x, tu(x)+(1-t)v(x), ... , tu^{(n)}(x)+(1-t)v^{(n)}(x) \right ) dt,  \\
& \vdots & \\
 c_{n-1}(x) & = & \int _0 ^1 \frac{\partial K}{\partial z_{n+1}} \left (x, tu(x)+(1-t)v(x), ... , tu^{(n)}(x)+(1-t)v^{(n)}(x) \right ) dt,  \\
c_{n}(x) & = & \int _0 ^1 \frac{\partial K}{\partial z_{n+2}} \left (x, tu(x)+(1-t)v(x), ... , tu^{(n)}(x)+(1-t)v^{(n)}(x) \right ) dt.
\end{eqnarray*}

\noindent
Let $w(x)=u(x)-v(x)$.  By (\ref{eq:functional}) we have
$$c_0 w+ c_1w' + \cdots + c_{n-1}w^{(n-1)} + c_{n} w^{(n)} \leq 0 .$$

\noindent
If $\frac{\partial K}{\partial z_{n+2}} >0$, then $c_n >0$, so
$$w^{(n)}(x)+ \frac{c_{n-1}}{c_n} w^{(n-1)} (x) + \cdots + \frac{c_{1}}{c_n}w'(x) + \frac{c_{0}}{c_n} w (x) \leq 0. $$

\vspace{.05in}
\noindent
The initial condition implies that
$$w(x_0) = 0, \quad w'(x_0) = 0, \quad ... \quad w^{(n-1)}(x_0) = 0.$$ 

\vspace{.05in}
\noindent
By Theorem \ref{thm:Hopf-equivalent}, there exists $\delta >0$ such that $w(x) \leq 0$ for $x \in (x_0, x_0 + \delta).$

\vspace{.1in}
\noindent
If $n$ is even, applying Theorem \ref{thm:Hopf-equivalent_right} and choosing a smaller $\delta$ if necessary, we know that $w(x) \leq 0$ for $x \in (x_0 - \delta, x_0).$

\vspace{.1in}
\noindent
If $n$ is odd, applying Theorem \ref{thm:Hopf-equivalent_right} and choosing a smaller $\delta$ if necessary, we know that $w(x) \geq 0$ for $x \in (x_0 - \delta, x_0).$

\vspace{.1in}
\noindent
Therefore,

\vspace{.05in}
\noindent
if $n$ is even, then $u(x) \leq v(x)$ for $x \in (x_0 - \delta, x_0 + \delta) ;$

\vspace{.05in}
\noindent
if $n$ is odd, then $u(x) \geq v(x)$ for $x \in (x_0 - \delta, x_0) $ and $u(x) \leq v(x)$ for $x \in (x_0 , x_0 + \delta) $.

\vspace{.05in}
\noindent
This completes the proof of Theorem \ref{thm:functional}.

\vspace{.2in}

\section{The Sequential Form of the Second Order Hopf's Lemma}
\label{sec:second_order}

\vspace{.05in}
Next, we will establish the higher order sequential versions of Hopf's lemma which are crucial in the proof of Theorem \ref{thm:functional}.  We first need to prove the following sequential Hopf's lemma in second order.

\vspace{.05in}

\begin{thm}
\label{thm:Hopf_k=2}
Let $u \in C^{2}((a,b)) \bigcap C^{1}([a,b))$ be a function which satisfies  
\begin{equation*}
\label{eq:main_k=2}
u''(x)+a_1(x)u'(x)+a_0(x)u(x) \leq 0  \hspace{.2in} \text{for } \,\, x\in (a,b)
\end{equation*}
 where  $|a_1(x)|$ and $|a_0(x)| $ are bounded by some constant $C>0$. Assume that $u$ satisfies 
\begin{equation*}
\label{eq:initial_k=2}
u(a)  =  0 ,
\end{equation*}
and Condition (\ref{eq:sequence}).
 
\noindent
Then $u'(a)>0$.  Furthermore, $u>0$ in a neighborhood of $a$.

\end{thm}

\vspace{.05in}
The classical second order Hopf's lemma requires that $u(x)>0$ for all $x$ greater than and sufficiently close to $a$, that is, $u(a)$ is a local minimum.  But here we only need the weaker assumption that $u$ is positive at a sequence of points approaching $a$, and we can show that then $u$ must be actually positive at all points near the boundary $a$.  In other words, $u(x)$ cannot oscillate around the $y$-axis as $x$ approaches $a$. 

\vspace{.05in}
In this section we present a proof of Theorem \ref{thm:Hopf_k=2} that relies on the following maximum principle on small intervals.  An alternative proof is given in the Appendix.

\begin{lemma}
\label{lem:small_interval}
Suppose $g \in C^{2}((a,b)) \bigcap C^{1}([a,b))$ satisfies  
$$L[g]=g''(x) +a_1(x)g'(x)+a_0(x)g(x) \leq 0 \hspace{.2in} \text{for } \,\, x\in (a,b), $$ where $|a_1(x)|,|a_0(x)|$ are bounded by some constant $C>0$. Then there exists a constant $\delta=\delta(C)>0$ such that on any interval $[c,d] \subseteq [a,b)$ with $|d-c|<\delta$, we have $g \geq 0$ provided $g(c) \geq 0$ and $g(d) \geq 0$. 

\end{lemma}

\vspace{.05in}

\pf
Without loss of generality we can assume $c=0$.  Define $$h(x)=e^{\gamma \delta}-e^{\gamma x} \hspace{.3in} \text{and} \hspace{.3in} w(x)=\frac{g(x)}{h(x)},$$ where $\gamma, \delta >0$ are to be chosen.  Then

\noindent 
\begin{eqnarray}
\label{eq:L[wh]}
L[g] & = &\frac{d^2}{dx^2}\left (w(x)h(x) \right )+a_1(x)\frac{d}{dx} \left ( w(x)h(x) \right )+a_0(x) \left ( w(x)h(x) \right ) \nonumber \\
& = & hw''+\left ( 2h'+a_1h \right )w'(x) + L[h]w(x).
\end{eqnarray}

\noindent
Suppose the minimum of $w$ is negative and achieved at some $x_0 \in (0, d)$.  Then $$w''(x_0) \geq 0, \qquad w'(x_0)=0, \qquad \text{and}\qquad w(x_0)<0. $$ 

\noindent
By definition $$h(x)>0 \qquad \text{if} \qquad 0\leq x \leq d < \delta.$$

\noindent
Direct computation shows that 
\begin{eqnarray*}
L[h] & = & e^{\gamma x} \cdot \left ( -\gamma^2 - a_1\gamma+a_0\left ( e^{\gamma \delta -\gamma x} -1 \right ) \right ) \\
& \leq & e^{\gamma x} \left ( -\gamma^2 + C \gamma + C (e^{\gamma \delta } -1) \right ) \hspace{.2in} \text{when} \hspace{.1in} 0\leq x \leq d < \delta .
\end{eqnarray*}

\noindent
We first choose $\gamma >0$ sufficient large so that $-\gamma ^2 + C \gamma + 2C <0$, then we choose $0 < \delta < \frac{\ln 3}{\gamma}$ so $0<e^{\gamma \delta } -1<2$.  Thus

$$L[h]  \leq  e^{\gamma x} \left ( -\gamma ^2 + C \gamma + 2 C \right) \,\, <  0 $$ when $ 0 \leq x \leq d < \delta $.

\vspace{.05in}

\noindent
Then by (\ref{eq:L[wh]}) it follows that  $L[g](x_0) >0$.  This contradiction proves that the minimum of $w$ on $[0,d]$ must be nonnegative, thus $g(x)\geq 0 $ on $[0,d]$ since $h(x)>0$.  

\qed

\vspace{.1in}
Next, we use Lemma \ref{lem:small_interval} to prove Theorem \ref{thm:Hopf_k=2}.  

\vspace{.05in}
\pf  Without loss of generality we can assume $a=0$.   

\vspace{.05in}
\noindent
Denote $$L[u]:=u''(x) +a_1(x)u'(x)+a_0(x)u(x).$$  

\noindent
Let $$g(x)=u(x)-\epsilon \left (e^{\lambda x}-1 \right ),$$ where $\epsilon >0$ will be chosen later.

\noindent
For $x \geq 0$ and $\lambda > 0$,
\begin{eqnarray*}
L[e^{\lambda x}-1] & = & \lambda ^2 e^{\lambda x} + a_1(x) \lambda e^{\lambda x} + a_0 (e^{\lambda x} -1) \\
& = & e^{\lambda x} \left ( \lambda ^2   + a_1 \lambda +  a_0 \left ( 1 - e^{- \lambda x} \right )\right )\\
& \geq & e^{\lambda x} \left ( \lambda ^2 - C \lambda - C \right ) \\
& > & 0
\end{eqnarray*}
when $\lambda$ is chosen to be sufficiently large.  Thus we know 
\begin{eqnarray*}
L[g] & = & L[u]-\epsilon L[e^{\lambda x}-1] \\
& < & 0 .
\end{eqnarray*}

\noindent
By definition $g(0)=0$.  Since the sequence $x_i \to 0$, we may choose an index $i_0$ such that $0<x_{i_0}< \delta$, where $\delta$ is chosen as in Lemma \ref{lem:small_interval}.  Because $u(x_{i_0})>0$, we can choose $$\epsilon = \frac{u(x_{i_0})}{e^{\lambda x_{i_0}}-1}>0$$ in the definition of $g(x)$.  Then we have $g(x_{i_0})=0 $.

\vspace{.05in}

\noindent
Now Lemma \ref{lem:small_interval} implies

$$g(x) \geq 0 \hspace{.3in} \text{on} \hspace {.3in} [0,x_{i_0}].$$

\noindent
The Taylor expansion of $g$ at $0$ gives $$g(x)=g'(0)x+O(x^2),$$ thus $g'(0) \geq 0$.  Consequently
$$
u'(0) \,\, = \,\, g'(0) + \epsilon \lambda \,\, > \,\, 0.
$$
 
\qed

\vspace{.1in}
Lemma \ref{lem:small_interval} shows that if $g$ is nonnegative at the two endpoints of a sufficiently small interval, then $g \geq 0$ in that interval.  For third and higher order differential inequalities, it no longer holds.  To see this, consider the sequence of functions $$g_i(x)= \left ( x - \frac{1}{i} \right )^2 - \frac{1}{i^2}.$$ Each function satisfies the differential equation $u_i^{(k)}=0$ for all $k=3,4,...$.  Although $g_i(0)=g_i(\frac{2}{i})=0$ and $\frac{2}{i} \to 0$, $g_i(x)$ is negative on $(0,\frac{2}{i})$.  

\vspace{.05in}
 
The classical maximum principle also fails in the higher order case.  For example, the function $u(x)=\sin x$ satisfies
\begin{eqnarray*}
u^{(3)}+u'+ 0 \cdot u & = & 0  \\
u^{(4)}+u''+ 0 \cdot u & = & 0 \\
& \vdots & 
\end{eqnarray*}
and $u(0)=u(2\pi)=0$, but $u \leq 0$ on $[\pi, 2 \pi]$.

\vspace{.05in}
 
Therefore, there exists a very interesting distinction between the Hopf's lemma and maximum principle in higher orders.  Although for the second order inequalities the Hopf's lemma can be used to prove the maximum principle, in the higher order case the maximum principle fails, but the Hopf's lemma still holds.

\vspace{.2in}

\section{The Higher Order Hopf's Lemmas}
\label{sec:higher_order}

\vspace{.05in}
Now we are ready to prove the higher order Hopf's Lemma, Theorems \ref{thm:Hopf} and \ref{thm:Hopf_right}.

\vspace{.1in}
\noindent
\textbf{Proof of Theorem \ref{thm:Hopf}}:

\vspace{.05in}
\noindent
We will employ a reduction of order technique and use mathematical induction. The case $n = 2$ is provided by Theorem \ref{thm:Hopf_k=2}.   Suppose the theorem is true for $n=k \geq 2$, we will show that it is also true for $n=k+1$, i.e. assume $u$ satisfies  
$$
u^{(k+1)}(x)+a_{k}(x)u^{k}(x)+\cdots +a_1(x)u'(x)+a_0(x)u(x) \leq 0  ,
$$
 where  $a_{k}(x), ... , a_1(x), a_0(x) $ are in $C([a,b))$,   
\begin{equation}
\label{eq:u_initial_induction}
u(a) \,\, = \,\, u'(a) \,\, = \,\, \cdots  \,\, = u^{(k-1)}(a) \,\, = \,\, 0 ,
\end{equation}
and Condition (\ref{eq:sequence}),
we need to show that $u^{(k)}(a)  >  0$.

\noindent
Let 
\begin{equation}
\label{eq:v_defn}
v:=fu+u',
\end{equation} 

\noindent
where $f$ is to be chosen.  We then have
\begin{eqnarray}
\label{eq:v_derivatives}
v' & = & f'u+fu'+u'' \nonumber \\
v'' & = & f''u+2f'u'+fu''+u^{(3)} \nonumber \\
v^{(3)} & = & f^{(3)}u + 3f''u' + 3f'u''+ fu^{(3)} +u^{(4)}\nonumber \\
& \vdots &    \\
v^{(k-2)} & = & f^{(k-2)} u+ \binom{k-2}{1}f^{(k-3)}u' + \cdots + \binom{k-2}{k-3} f'u^{(k-3)}+ fu^{(k-2)} +u^{(k-1)}  \nonumber \\
v^{(k-1)} & = & f^{(k-1)} u+ \binom{k-1}{1}f^{(k-2)}u' + \cdots + \binom{k-1}{k-2} f'u^{(k-2)}+ fu^{(k-1)} +u^{(k)}  \nonumber \\
v^{(k)} & = & f^{(k)} u+ \binom{k}{1}f^{(k-1)}u' + \cdots + \binom{k}{k-1} f'u^{(k-1)}+ fu^{(k)} +u^{(k+1)} . \nonumber 
\end{eqnarray}

\vspace{.1in}
\noindent
We would like to choose appropriate functions $b_0(x), b_1(x), ... , b_{k-1}(x) \in C([a,b))$, such that 
\begin{eqnarray}
\label{eq:Lv_Lu}
& & u^{(k+1)}(x) + a_{k}(x)u^{k}(x)+\cdots +a_1(x)u'(x)+a_0(x)u(x) \nonumber \\
& = &  v^{(k)}(x) + b_{k-1}(x)v^{k-1}(x)+\cdots +b_1(x)v'(x)+b_0(x)v(x) .
\end{eqnarray}

\vspace{.05in}
\noindent
Because of (\ref{eq:v_defn}) and (\ref{eq:v_derivatives}), the right hand side of (\ref{eq:Lv_Lu}) becomes
\begin{eqnarray*}
& & \hspace{.5in} f^{(k)} u+ \binom{k}{1}f^{(k-1)}u' + \cdots + \binom{k}{k-1} f'u^{(k-1)}+ fu^{(k)} +u^{(k+1)} \\
& & +b_{k-1} \left [  f^{(k-1)} u+ \binom{k-1}{1}f^{(k-2)}u' + \cdots + \binom{k-1}{k-2} f'u^{(k-2)}+ fu^{(k-1)} +u^{(k)} \right ] \\
& & + b_{k-2} \left [  f^{(k-2)} u+ \binom{k-2}{1}f^{(k-3)}u' + \cdots + \binom{k-2}{k-3} f'u^{(k-3)}+ fu^{(k-2)} +u^{(k-1)} \right ] \\
& & + \cdots + b_2 \left ( f''u+2f'u'+fu''+u^{(3)} \right ) + b_1 \left ( f'u+fu'+u'' \right ) + b_0(fu+u'),
\end{eqnarray*}

\vspace{.1in}

\noindent
which is equal to
\begin{eqnarray*}
& & u^{(k+1)} + \left ( f+b_{k-1} \right )u^{(k)} + \left [ \binom{k}{k-1}f'+b_{k-1}f+b_{k-2}   \right ] u^{(k-1)} \\
& & + \left [ \binom{k}{k-2}f''+b_{k-1}\binom{k-1}{k-2} f' + b_{k-2}f +  b_{k-3} \right ] u^{(k-2)} + \cdots \\
& & + \left [ \binom{k}{2}f^{(k-2)} + b_{k-1}\binom{k-1}{2} f^{(k-3)}+ b_{k-2}\binom{k-2}{2} f^{(k-4)}+ \cdots +  b_{2}f +  b_1 \right ] u'' \\
& & + \left [ \binom{k}{1}f^{(k-1)} + b_{k-1}\binom{k-1}{1} f^{(k-2)}+ b_{k-2}\binom{k-2}{1} f^{(k-3)}+ \cdots +  b_{1}f +  b_0 \right ] u' \\
& & + \left ( f^{(k)} + b_{k-1} f^{(k-1)}+ b_{k-2} f^{(k-2)}+ \cdots +  b_{1}f' +  b_0 f \right ) u.
\end{eqnarray*}

\vspace{.1in}
\noindent
In light of (\ref{eq:Lv_Lu}), we want to choose $b_0(x), b_1(x), ... , b_{k-1}(x)$ such that
\begin{eqnarray}
\label{eq:a_expressions}
a_{k} & = & f+b_{k-1}  \nonumber \\
a_{k-1} & = & \binom{k}{k-1}f'+b_{k-1}f+b_{k-2}   \nonumber \\
a_{k-2} & = & \binom{k}{k-2}f''+b_{k-1}\binom{k-1}{k-2} f' + b_{k-2}f +  b_{k-3}  \nonumber  \\
& \vdots & \\
a_{2} & = &\binom{k}{2}f^{(k-2)} + b_{k-1}\binom{k-1}{2} f^{(k-3)}+ b_{k-2}\binom{k-2}{2} f^{(k-4)}+ \cdots +  b_{2}f +  b_1  \nonumber \\
a_{1} & = & \binom{k}{1}f^{(k-1)} + b_{k-1}\binom{k-1}{1} f^{(k-2)}+ b_{k-2}\binom{k-2}{1} f^{(k-3)}+ \cdots +  b_{1}f +  b_0  \nonumber \\
a_{0} & = & f^{(k)} + b_{k-1} f^{(k-1)}+ b_{k-2} f^{(k-2)}+ \cdots +  b_{1}f' +  b_0 f.  \nonumber 
\end{eqnarray} 

\vspace{.05in}
\noindent
Solving for $b_{k-1}, ... , b_1, b_0$ from the first $k$ equations, we obtain
\begin{eqnarray}
\label{eq:b_k-1_to_0}
b_{k-1} & = & a_{k} - f \nonumber \\
b_{k-2} & = & a_{k-1} - \binom{k}{k-1}f' - b_{k-1}f \nonumber\\
b_{k-3} & = & a_{k-2} - \binom{k}{k-2}f'' - b_{k-1} \binom{k-1}{k-2}f'- b_{k-2} f \nonumber \\
& \vdots &   \\
b_1 & = & a_2 -  \binom{k}{2}f^{(k-2)} - b_{k-1}\binom{k-1}{2} f^{(k-3)} - b_{k-2}\binom{k-2}{2} f^{(k-4)} - \cdots -  b_{2}f   \nonumber \\
b_0 & = & a_{1} - \binom{k}{1}f^{(k-1)} - b_{k-1}\binom{k-1}{1} f^{(k-2)} - b_{k-2}\binom{k-2}{1} f^{(k-3)} - \cdots -  b_{1}f . \nonumber 
\end{eqnarray}

\noindent
If the first equation in (\ref{eq:b_k-1_to_0}) is substituted into the second equation, $b_{k-2}$ can be expressed as  $a_{k-1} - \binom{k}{k-1}f' - a_k f + f^2$, which is a polynomial in $f$ and $f'$ with coefficients comprised of $a_{k}$, $a_{k-1} $ and universal constants.  Similarly $b_{k-3}$, ... , $b_1$, $b_0$ all can be expressed as polynomials in $f$ and its derivatives, with the coefficients given by $a_0(x)$, ... , $a_k(x)$ and universal constants. 

\vspace{.05in}
\noindent
Thus we can write
\begin{eqnarray}
\label{eq:b_expressions}
b_{k-1} & = &  P_{k-1} \left (a_{k}, f \right ) \nonumber \\
b_{k-2} & = &  P_{k-2} \left (a_{k}, a_{k-1}, f, f' \right )  \nonumber \\
b_{k-3} & = &  P_{k-3} \left (a_{k}, a_{k-1}, a_{k-2}, f, f', f'' \right )  \nonumber \\
& \vdots & \\
b_1 & = & P_{1} \left (a_{k}, a_{k-1}, ... , a_{2}, f, f', ... , f^{(k-2)} \right ) \nonumber \\
b_0 & = & P_{0} \left (a_{k}, a_{k-1}, ... , a_{1}, f, f', ... , f^{(k-1)} \right ). \nonumber 
\end{eqnarray}

\noindent
Here $P_{k-1}$, $P_{k-2}$, ... , $P_{1}$, $P_{0}$ are polynomials in $f$ and its derivatives, and their coefficients depend on the continuous functions $a_{k}(x)$, $a_{k-1}(x)$, ... , $ a_{1}(x)$.  

\vspace{.05in}
\noindent
Then we substitute (\ref{eq:b_expressions}) into the last equation in (\ref{eq:a_expressions}), so the function $f$ must satisfy the $k$-th order ODE
\begin{equation}
\label{eq:f}
f^{(k)}+ P_{k-1}f^{(k-1)} + P_{k-2}f^{(k-2)} + \cdots + P_1 f'+P_0 f = a_0 .
\end{equation} 

\vspace{.05in}
\noindent
Under the initial condition $f(a)=1$, Equation (\ref{eq:f}) has a solution $f \in C^k([a,a+\epsilon))$ for some $\epsilon>0$.  With this choice of $f$, (\ref{eq:Lv_Lu}) holds, so we know that 
$$v^{(k)}(x) + b_{k-1}(x)v^{k-1}(x)+\cdots +b_1(x)v'(x)+b_0(x)v(x) \leq 0 . $$ 

\vspace{.05in}
\noindent
Definition (\ref{eq:b_expressions}) implies that the coefficient functions $b_{k-1}(x)$, ... $b_1(x)$, $b_0(x)$ are all continuous. 

\vspace{.05in}
\noindent
Since $f(a)=1$ and 
\begin{equation*}
u(a)=u'(a)=...=u^{k-1}(a)=0,
\end{equation*} 
from (\ref{eq:v_derivatives}) we know that $$v(a)=v'(a)=...=v^{(k-2)}(a)=0.$$
 
\vspace{.05in}
\noindent
Because there exists a sequence $x_i \to a$ with $u(x_i)>0$ and $u(a)=0$, we can choose a sequence $\tilde{x}_i \to a$ such that $u(\tilde{x}_i)>0$ and $u'(\tilde{x}_i)>0$.  Since $f(a)=1$, when $i$ is sufficiently large we have $f(\tilde{x}_i)>0$.  Therefore $$v(\tilde{x}_i)=f(\tilde{x}_i)u(\tilde{x}_i)+u'(\tilde{x}_i) >0.$$

\vspace{.05in}
\noindent
Thus by the inductive hypothesis we know $$v^{(k-1)}(a)>0.$$

\noindent
Then the second last equation in (\ref {eq:v_derivatives}) and the initial conditions (\ref {eq:u_initial_induction}) implies $$u^{(k)}(a)>0.$$  The proof of Theorem \ref{thm:Hopf} is now completed by mathematical induction.

\qed

\vspace{.1in}
\noindent
\textbf{Proof of Theorem \ref{thm:Hopf_right}}:

\begin{enumerate} [(i)]

\item 

If $n$ is even, define $$\hat{u}(x):=u(2b-x).$$  Then $\hat{u} \in C^n((b, 2b-a)) \bigcap C^{n-1}([b, 2b-a))$ and 
\begin{eqnarray*}
\hat{u}'(x) & = & - u'(2b-x) \\
\hat{u}''(x) & =  &  u''(2b-x) \\
& \vdots & \\
\hat{u}^{(n-1)}(x) & = & (-1)^{n-1}u^{(n-1)}(2b-x) \\
& = & - u^{(n-1)}(2b-x) \\
\hat{u}^{(n)}(x) & = & (-1)^{n}u^{(n)}(2b-x)\\
& = & u^{(n)}(2b-x), 
\end{eqnarray*}
\noindent
and $\hat{u}$ satisfies
\begin{equation*}
\hat{u}^{(n)}(x)-a_{n-1}(2b-x)\hat{u}^{n-1}(x)+\cdots - a_1(2b-x)\hat{u}'(x) + a_0(2b-x) \hat{u}(x) \leq 0 ,
\end{equation*}

\noindent
where the functions $ a_0(2b-x), -a_1(2b-x), ... , -a_{n-1}(2b-x)$ are in $C([b,2b-a))$.  

\vspace{.05in}
\noindent
The initial conditions (\ref {eq:initial_right}) imply that 
\begin{equation*}
\hat{u}(b) \,\, = \,\, \hat{u}'(b) \,\, = \,\, \cdots \,\, = \,\, \hat{u}^{(n-2)}(b) \,\, = \,\, 0.
\end{equation*}

\noindent
By (\ref{eq:sequence_even}), there exists a sequence $\{ 2b-x_i \}$, such that $b<2b-x_i <2b-a$, $2b-x_i\to b$, and $\hat{u} (2b-x_i) =  u(x_i) > 0.$

\vspace{.05in}
\noindent
Then by Theorem \ref{thm:Hopf},  $\hat{u}^{(n-1)}(b)>0$ and $\hat{u}>0$ in a neighborhood of $b$.  Therefore,  we have $u^{(n-1)}(b) < 0$ and $u > 0$ in a neighborhood of $b$.

\vspace{.05in}

\item

If $n$ is odd, define $$\tilde{u}(x):=-u(2b-x).$$  Then $\tilde{u} \in C^n((b, 2b-a)) \bigcap C^{n-1}([b, 2b-a)$ and 
\begin{eqnarray*}
\tilde{u}'(x) & = & u'(2b-x) \\
\tilde{u}''(x) & =  &  - u''(2b-x) \\
& \vdots & \\
\tilde{u}^{(n-1)}(x) & = & (-1)^{n}u^{(n-1)}(2b-x) \\
& = & - u^{(n-1)}(2b-x) \\
\tilde{u}^{(n)}(x) & = & (-1)^{n+1}u^{(n)}(2b-x)\\
& = & u^{(n)}(2b-x), 
\end{eqnarray*}

\noindent
and $\tilde{u}$ satisfies
\begin{equation*}
\tilde{u}^{(n)}(x)-a_{n-1}(2b-x)\tilde{u}^{n-1}(x)+\cdots + a_1(2b-x)\tilde{u}'(x) - a_0(2b-x) \tilde{u}(x) \leq 0 ,
\end{equation*}

\noindent
where the functions $- a_0(2b-x), a_1(2b-x), ... , -a_{n-1}(2b-x)$ are in $C([b,2b-a))$.  

\vspace{.05in}
\noindent
The initial conditions (\ref {eq:initial_right}) imply that 
\begin{equation*}
\tilde{u}(b) \,\, = \,\, \tilde{u}'(b) \,\, = \,\, \cdots \,\, = \,\, \tilde{u}^{(n-2)}(b) \,\, = \,\, 0.
\end{equation*}
\noindent
By (\ref{eq:sequence_odd}), there exists a sequence $\{ 2b-x_i \}$, such that $b<2b-x_i <2b-a$, $2b-x_i\to b$, and $\tilde{u} (2b-x_i) = - u(x_i) > 0$.

\vspace{.05in}
\noindent
Then by Theorem \ref{thm:Hopf},  $\tilde{u}^{(n-1)}(b)>0$ and $\tilde{u}>0$ in a neighborhood of $b$.  Therefore, we have $u^{(n-1)}(b) < 0$ and $u < 0$ in a neighborhood of $b$.

\end{enumerate}

\qed

\vspace{.2in}
\section{Some Comments on the Proofs of Higher Order Hopf's Lemma}
\label{sec:comments}

\vspace{.1in}
The proof of Theorem \ref{thm:Hopf} shows that it is necessary to first obtain the sequential form of the second order Hopf's lemma (Theorem \ref{thm:Hopf_k=2}), as we only know the sign of the function $v$ at a sequence of points after the reduction process, so the classical Hopf's lemma no longer applies. 

\vspace{.1in}
It is worth pointing out that the conditions (\ref {eq:sequence}), (\ref {eq:sequence_even}), and (\ref {eq:sequence_odd}) are sharp in the sense that if they are not satisfied, then the $(n-1)$-th derivative may vanish at the endpoints.

\vspace{.1in}
\noindent
\textbf{Example}:
For any $ 0<\alpha <1$ and $n \geq 3,$ define 
$$u = \begin{cases}(-1)^{n-1}\lambda_n (-x)^{\frac{ n}{1-\alpha}},  & x <0 \\
     \qquad -  ~\lambda_n x^{\frac{n}{1-\alpha}},  & x \geq 0
      \end{cases} $$
where $$ \lambda_n = \left[ (\beta+n) \cdots(\beta+1)\right]^{\frac{1}{\alpha-1}}, \qquad \text{and} \qquad \beta=\frac{ n}{1-\alpha}-n=\frac{n\alpha }{1-\alpha}.$$

\noindent
Direct computation shows that
\begin{equation}
\label{eq:Holder}
u^{(n)}(x) = -|u(x)|^\alpha, \qquad x\in(-\infty,\infty) .
\end{equation}
Therefore $u$ satisfies the differential inequality $$u^{(n)} \leq 0.$$

\noindent
To simplify the expressions let us choose $\alpha = \frac{1}{2}$, then
$$u = 
\begin{cases} (-1)^{n-1} \left(\frac{n!}{(2n)!}\right)^2 (-x)^{2n},  & x <0 \\
\qquad -  \left(\frac{n!}{(2n)!}\right)^2 x^{2n},  & x \geq 0.
\end{cases} 
$$
By definition $$ u(0) = u'(0) = \cdots = u^{(n-2)}(0)=0$$ and also $$u^{(n-1)}(0)=0 .$$  

\vspace{.05in}
\noindent  
Note that $u < 0$ on $(0,1)$, so Condition (\ref {eq:sequence}) is not satisfied on $(0,1)$.

\vspace{.05in}
\noindent    
If $n$ is even, $u < 0$ on $(-1,0)$, so Condition (\ref {eq:sequence_even}) is not satisfied on $(-1,0)$.

\vspace{.05in}
\noindent    
If $n$ is odd, $u > 0$ on $(-1,0)$, so Condition (\ref {eq:sequence_odd}) is not satisfied on $(-1,0)$.
 
\qed

\vspace{.1in}
Theorems \ref{thm:Hopf} and \ref{thm:Hopf_right} need to assume that the coefficient functions $a_0(x)$,..., $a_{n-1}(x)$ are continuous, while in Theorem \ref{thm:Hopf_k=2} they only need to be bounded.  The continuity condition is assumed when $n \geq 3$ to ensure that Equation (\ref {eq:f}) possesses a solution $f$. It would be interesting to know whether this is merely a limitation of the technique used in the proof or this reflects an inherent difference between the second and higher order cases.  When $n=3$, the continuity requirement can be replaced by boundedness, if we assume an additional assumption that $u$ be non-negative at all points near $a$.

\begin{thm}
\label{thm:Hopf_k=3}
Let $u\in\mathcal C^3((a,b)) \bigcap \mathcal C^2([a,b))$ be a function that satisfies
\begin{equation*} 
u^{(3)}(x) + a_2(x) u''(x) + a_1(x) u'(x) + a_0(x) u(x) \leq 0 \hspace{.2in} \text{for } \,\, x\in (a,b),
\end{equation*}
where $|a_0(x)|, |a_1(x), |a_2(x)|\leq C$ for some constant $C>0$. Suppose $u(a)=u'(a)=0$, $u(x) \geq 0$ for all $x$ in a small neighborhood of $a$, and there exists a sequence $\{x_i \} \subset (a,b)$ such that $x_i \to a$ and $u(x_i)>0$.  
Then   
$ u''(a)>0. $
\end{thm} 

\pf If $a_0(x) \geq 0$ for $x$ in a small neighborhood of $a$, then since $u(x) \geq 0$ near $a$, we have
$$ v''(x) + a_2(x) v'(x) + a_1(x) v(x) \leq 0 \qquad \text{for}   \quad x\in(a, a+\epsilon) \subset (a,b) , $$
where $$v(x) = u'(x)  \qquad \text{and} \qquad v(a) = u'(a) = 0.$$  

\vspace{.05in}
\noindent
Suppose $v(x)\leq 0$ for all $x$ near $a$, then $v(x)=u'(x)$ and $u(a)=0$ imply $u(x)\leq 0$ on $(a, a + \epsilon)$, contradicting the assumption that $x_i \to a$ and $u(x_i)>0$.  Therefore there exists a sequence $\tilde{x}_i \to a$ such that $v(\tilde{x}_i) > 0$.  By Theorem \ref{thm:Hopf_k=2}, we then have $v'(a)=u''(a)>0$.

\vspace{.1in}
\noindent 
For general $a_0(x)$, let 
$$ m(x) = e^{\theta \eta}-e^{- \theta (x-a)} \quad \text{for} \quad a \leq x \leq a+\eta < b. $$
For each $\theta>0$ we may choose $\eta$ such that 
\begin{equation}
\label{eq: eta}
e^{2\theta\eta}-1<\theta(b-a) .
\end{equation}
Then since $e^{\theta(\eta+x - a)} \leq e^{2\theta \eta}$, we have
$$ e^{\theta(\eta+x - a)} = 1 + h(x), \qquad \text{where} \quad 0< h(x) < \theta(b-a) \,\, \text{for all} \,\ x\in (a, a+\eta). $$
Because $|a_2(x)|, |a_1(x)|, |a_0(x)|\leq C$,   for $a<x<a+\eta$
\begin{eqnarray*}
L[m] & := & m^{(3)}(x) + a_2(x) m''(x) + a_1(x) m'(x) + a_0(x)  m(x) \\
&= &
\left( \theta^3 - a_2(x)\theta^2 + a_1(x) \theta + a_0(x) (e^{\theta (\eta+ x-a)} -1)\right)e^{-\theta (x-a)}  \\
&= &
\left( \theta^3 - a_2(x)\theta^2 + a_1(x) \theta + a_0(x) h(x) \right)e^{-\theta (x-a)}  \\
& \geq &
\left( \theta^3 - |a_2(x)|\theta^2 - |a_1(x)| \theta - |a_0(x)| \theta(b-a)  \right)e^{-\theta (x-a)} \\
& \geq &
\left( \theta^3 -C\theta^2 - (1+b-a)C \theta\right)e^{-\theta (x-a)} . 
\end{eqnarray*} 

\vspace{.05in}
\noindent
We can choose $\theta$ to be sufficiently large such that $$ \theta^3 -C\theta^2 - (1+b-a)C \theta >0 .$$  With this  $\theta$, choose $\eta$ as above to satisfy (\ref {eq: eta}).  Then we have $$L[m] >0.$$

\vspace{.05in}
\noindent
For $x\in[a, a + \eta]$,  $m(x)>0$ by definition, so we may define 
$$z(x) =\frac{u(x)}{m(x)}. $$

\noindent
Applying the differential operator $L$ to $u(x) = m(x) z(x)$, 
\begin{eqnarray*}
L[u] &= & (m(x) z(x))^{(3)} + a_2(x)  (m(x) z(x))'' + a_1(x)  (m(x) z(x))' + a_0(x)  (m(x) z(x)) \\
&= & m(x) z^{(3)}(x) + [3 m'(x) + a_2(x) m(x)] z''(x) \\
& & + [3 m''(x) + 2 a_2(x) m'(x)+ a_1(x)m(x)] z'(x) + L[m]  z(x) .
\end{eqnarray*}
Since $L[u]\leq 0$ and  $m(x) >0$, we have 
\begin{equation}
\label{eq:z}
 z^{(3)}(x) +a_2^*(x) z''(x) + a_1^*(x) z'(x) + a_0^*(x) z(x)  \leq 0  ,
\end{equation}
where
\begin{eqnarray*}
a_2^*(x)&  = & \frac{3 m'(x)}{m(x)} + a_2(x), \\
a_1^*(x) & = & \frac{3 m''(x) + 2 a_2(x) m'(x)}{m(x)} + a_1(x), \\
a_0^*(x) & = & \frac{L[m]}{m(x)} 
\end{eqnarray*}

\noindent
For fixed $\theta$, $\frac{ m'}{m}, \frac{m''}{m}$ and $\frac{L[m]}{m}$ are all bounded when $x\in(a, a + \eta]\subset (a,b)$, so there exists $C_1>0$ such that
$$ |a_2^*(x)|, \,\, |a_1^*(x)|, \,\, |a_0^*(x)| \leq C_1. $$
Since $u'(a)=u(a)=0$ and $m'(a) = \theta \neq 0$, we have
$$ \lim_{x\to a+} z'(x) =  \lim_{x\to a+} \frac{u'(x)m(x) - u(x) m'(x)}{m^2(x)}  = \frac{u'(a)m(a) - u(a) m'(a)}{m^2(a)} = 0.$$
The function $z(x)\in \mathcal C^2([a, a + \eta])$ satisfies
$$ z'(a) = z(a)=0, \qquad z(x_i)=\frac{u(x_i)}{m(x_i)}>0, \qquad\text{and} \quad z(x) \geq 0 \quad \text{for} \quad x\in (a, a + \eta]. $$
Recall that $m>0$ and $L[m]>0$, so $a_0^*(x) >0$ on $[a, a + \eta]$. Then by (\ref{eq:z}) and the discussion at the beginning of this proof we conclude that
$$ z''(a) > 0. $$
Consequently,
\begin{eqnarray*}
 u''(a) & = & m''(a) z(a) + 2 m'(a) z'(a) + m(a) z''(a) \\
 & = & m(a) z''(a) \\
 & > & 0. 
 \end{eqnarray*}
 This completes the proof.

\qed

\vspace{.05in}
It is natural to ask if Theorems \ref{thm:Hopf} and \ref{thm:Hopf_right} can be generalized to include two or more variables.  Generally speaking the answer is no.  Even the second order sequential Hopf's lemma fails with two variables.  For example, the function $u(x,y)=xy$ satisfies $\Delta  u =0$.  Although $u(0,0)=0$ and we can find a sequence of points $(x_i,y_i) \to (0,0)$ with $u(x_i,y_i)>0$, all directional derivatives of $u$ vanish at $(0,0)$ because $\nabla u (0,0)= (0,0).$

\vspace{.05in}
It also seems to be difficult to correctly formulate a multi-variable version of a higher order Hopf's lemma.  When $n$ is odd, Conditions (\ref {eq:sequence}) and (\ref {eq:sequence_odd}) require $u(x_i)$ to assume different sign at the two endpoints, and $u^{(n-1)}(a)$ and $u^{(n-1)}(b)$ have opposite sign in Theorems \ref{thm:Hopf} and \ref{thm:Hopf_right}.  

\vspace{.05in}
\noindent
This ``boundary effect'' is not an issue when $n=2$ because it is an even number and $u'(b)= - D_{\eta} u (b)$, where $\eta$ denotes the direction pointing toward the center of the interval.  Therefore, Theorems \ref{thm:Hopf} and \ref{thm:Hopf_right} and be combined to state that $D_{\eta} u >0$ on the boundary of the interval $(a,b)$.  When $n$ is odd, however, we will not be able to unify the two derivatives at the two endpoints.  In the multi-variable case, the boundary will be even more complicated, so it appears to be difficult to formulate a clear and unified expression for the derivatives like the one in the classical Hopf's lemma.

\vspace{.1in}

\vspace{.2in}
\section{Applications of Higher Order Hopf's Lemmas}
\label{sec:applications}

\vspace{.1in}
In this section we will give some additional applications of the higher order Hopf's lemmas.

\vspace{.1in}
Applying Theorem \ref{thm:Hopf-equivalent} to both functions $u$ and $-u$ gives a new proof of the standard uniqueness theorem of linear ODEs:

\vspace{.05in}

\begin{cor}
\label{thm:Hopf-corollary}
Let $u \in  C^{n}((a,b)) \bigcap C^{n-1}([a,b)) $ be a function which satisfies  
\begin{equation*}
\label{eq:main_equal_0}
u^{(n)}(x)+a_{n-1}(x)u^{n-1}(x)+\cdots +a_1(x)u'(x)+a_0(x)u(x) = 0  \hspace{.2in} \text{for } \,\, x\in (a,b)
\end{equation*}
 where $n\geq 2$ is a positive integer and $a_{n-1}(x), ... , a_1(x), a_0(x) $ are in $C([a,b))$. Assume that $u$ satisfies 
\begin{equation*}
\label{eq:initial_cor}
u(a) \,\, = \,\, u'(a) \,\, = \,\, \cdots  \,\, = u^{(n-2)}(a) \,\, = u^{(n-1)}(a) \,\, =\,\, 0 .
\end{equation*}

\noindent
Then $u \equiv 0$.

\vspace{.05in}
\noindent
When $n=2$, it suffices to assume that $a_1(x)$ and $a_0(x) $ are bounded functions.

\end{cor}

\vspace{.05in}

Another immediate consequence of Theorem \ref{thm:Hopf} is a unique continuation theorem.

\begin{cor}
\label{thm:unique_continuation}
Suppose  $u \in C^{\infty}([a,b))$ satisfies
$$u^{(n)}(x)+a_{n-1}(x)u^{n-1}(x)+\cdots +a_1(x)u'(x)+a_0(x)u(x) \leq 0, $$
where $n\geq 2$ is a positive integer and $a_{n-1}(x), ... , a_1(x), a_0(x) $ are in $C([a,b))$.  If Condition (\ref{eq:sequence}) holds, then it cannot be true that $u^{(k)}(a)=0$ for all $ k=0, 1, ....$

\vspace{.05in}
\noindent
When $n=2$, it suffices to assume that $a_1(x)$ and $a_0(x) $ are bounded functions.

\end{cor}

\vspace{.05in}

When $u$ is in $C^{\infty}([a,b))$, Theorem \ref{thm:Hopf} also follows from Corollary \ref{thm:unique_continuation}, hence the two results are equivalent.   Here is the proof.

\vspace{.1in}
\noindent 
\pf  Assume Corollary \ref{thm:unique_continuation} holds and $u \in C^{\infty}([a,b))$ satisfies (\ref{eq:main}), (\ref{eq:initial}), and (\ref{eq:sequence}), we need to show that $u^{(n-1)}(a)>0$.

\vspace{.05in}
\noindent
Condition (\ref{eq:sequence}) and the $(n-1)$-th degree Taylor's expansion of $u$ near $a$ implies that $u^{(n-1)}(a) \geq 0$.  

\vspace{.05in}
\noindent
Suppose $u^{(n-1)}(a)=0$.  

\vspace{.05in}
\noindent
Then by the $n$-th degree Taylor's expansion of $u$ near $a$ we have $u^{(n)}(a) \geq 0 $.  On the other hand, (\ref{eq:main}) and (\ref{eq:initial}) imply $u^{(n)}(a) \leq 0 $.  Therefore $u^{(n)}(a) = 0 $.  Again the $(n+1)$-th degree Taylor's expansion of $u$ near $a$ implies that $u^{(n+1)}(a) \geq 0$.  

\vspace{.05in}
\noindent
If $u^{(n+1)}(a) > 0$, then for $x$ close to $a$,
\begin{eqnarray*}
u(x) & = & \frac{u^{(n+1)}(a)}{(n+1)!}(x-a)^{n+1} + O \left (  (x-a)^{n+2} \right ) \\
u'(x) & = &  \frac{u^{(n+1)}(a)}{n!}(x-a)^{n} + O \left (  (x-a)^{n+1} \right ) \\
& \vdots & \\
u^{(n-1)}(x) & = & \frac{u^{(n+1)}(a)}{2!}(x-a)^2 + O \left (  (x-a)^3 \right ) \\
u^{(n)}(x) & = & \frac{u^{(n+1)}(a)}{1!}(x-a) + O \left (  (x-a)^2 \right ) .
\end{eqnarray*}

\noindent
Therefore, since $a_0(x), ... , a_{n-1}(x)$ are bounded,
\begin{eqnarray*}
& & u^{(n)}(x)+a_{n-1}(x)u^{n-1}(x)+\cdots +a_1(x)u'(x)+a_0(x)u(x) \\
& = & u^{(n+1)}(a)  \Big [ (x-a) + \frac{a_{n-1}(x)}{2!} (x-a)^2 + \cdots +  \frac{a_2(x)}{(n-1)!} (x-a)^{n-1} \\
& & + \frac{a_1(x)}{n!} (x-a)^{n}  + \frac{a_0(x)}{(n+1)!} (x-a)^{n+1}   \Big ] + O \left (  (x-a)^2 \right )  \\
& > & 0,
\end{eqnarray*}
provided that $x-a>0$ is sufficiently small.  This contradicts (\ref {eq:main}).  Hence $u^{(n+1)}(a) = 0$.

\vspace{.05in}
\noindent
Next we can show by similar argument that $u^{(n+2)}(a)=0$, then $u^{(n+3)}(a)=0$, ....  So $u^{(k)}(a)=0$ for all $k=0,1,2...$.  But this contradicts Corollary \ref {thm:unique_continuation}.  Therefore we must have $u^{(n-1)}(a)>0$, and Theorem \ref{thm:Hopf} holds.

\qed

\vspace{.1in}

The last application is about the boundary behavior of solutions to a type of nonlinear ODEs.  A similar ``boundary estimate'' concerning solutions of boundary-value problem for a semilinear Poisson PDE was given in \cite{Evans}. 

\begin{thm}
\label{lem:boundary}  
Let $u \in C^n([a,b])$ satisfy 
\begin{equation}
\label{eq:boundary_lemma}
 u^{(n)}(x) =  f(u, u',..., u^{(n-1)})  \quad \text{in $[a,b]$ },
\end{equation}
where $f(z_1, ... , z_{n}): \mathbf{R} ^{n} \rightarrow  \mathbf{R}$ is Lipschitz continuous in all variables.

\begin{enumerate} [(i)]

\item  Assume $u(a) = u'(a) = \cdots = u^{(n-2)}(a) = 0 $ and $u>0$ in a neighborhood of $a$.  Then either
$$u^{(n-1)}(a)>0$$ or $$u^{(n-1)}(a)=0, \quad u^{(n)}(a)>0.$$

\noindent
In either case, $u$ is strictly increasing near $a$.

\vspace{.05in}
\item  Assume $u(b) = u'(b) = \cdots = u^{(n-2)}(b) = 0 $ and $u>0$ in a neighborhood of $b$.  Then either
$$(-1)^{n-1}u^{(n-1)}(b) > 0$$ or $$u^{(n-1)}(b)=0, \quad (-1)^n u^{(n)}(b) > 0. $$

\noindent
In either case, $u$ is strictly decreasing near $b$.

\end{enumerate}

\end{thm}

\pf
 
\begin{enumerate} [(i)]

\item Assume $u(a) = u'(a) = \cdots = u^{(n-2)}(a) = 0 $ and $u>0$ in a neighborhood of $a$. 

\vspace{.1in}

\noindent
\textsl{Case 1}:  \,\,  $f(0,0,..., 0) \leq 0$.  

\vspace{.05in}
\noindent
Since $f$ is Lipschitz, it is differentiable almost everywhere.  Then from (\ref {eq:boundary_lemma}) we have
\begin{eqnarray*}
f(0,0,..., 0)  & = & \left ( f(0,0,..., 0) - f(u, u',..., u^{(n-1)}) \right ) + u^{(n)}(x)  \\
 & = & -\left ( \int _0 ^1 \frac{\partial f}{\partial z_{1}}(tu, tu',..., tu^{(n-1)}) \,\, dt \right)  u \\
        &  & - \left ( \int _0 ^1 \frac{\partial f}{\partial z_{2}}(tu, tu',..., tu^{(n-1)}) \,\, dt \right )  u' \\
        &  & - \cdots  \\
& & - \left ( \int _0 ^1 \frac{\partial f}{\partial z_{n}}(tu, tu',..., tu^{(n-1)}) \,\, dt \right )  u^{(n-1)} + u^{(n)}(x)  .
\end{eqnarray*}

\noindent
Hence $u$ satisfies $$u^{(n)}+ a_{n-1}  u^{(n-1)} + \cdots  + a_{1}  u' +  a_{0}   u  = f(0, 0, ... , 0) \leq 0 , $$

\noindent
where  
\begin{eqnarray*}
a_{n-1}(x) & = & - \displaystyle \int _0 ^1 \frac{\partial f}{\partial z_{n}}(tu, tu',..., tu^{(n-1)}) \,\, dt, \\
& \vdots & \\
a_0(x) & = & - \displaystyle \int _0 ^1 \frac{\partial f}{\partial z_{1}}(tu, tu',..., tu^{(n-1)}) \,\, dt.
\end{eqnarray*}  

\noindent
By Theorem \ref{thm:Hopf}, we have $$u^{(n-1)}(a)>0.$$

\vspace{.05in}
\noindent
The $(n-3)$-th degree Taylor expansion of $u'(x)$ near $a$ gives
$$u'(x)=\frac{1}{(n-2)!}u^{(n-1)}(\theta)(x+1)^{n-2}  \hspace{.2in} \text{for some} \hspace{.1in} a<\theta<x .$$

\noindent
When $x$ is sufficiently close to $a$, $u^{(n-1)}(\theta)>0$.  Thus $u'(x)>0$ and $u(x)$ is strictly increasing.

\vspace{.1in}
\noindent
\textsl{Case 2}: \,\, $f(0,0,..., 0) > 0$.  

\vspace{.05in}
\noindent
Since $u(x)>0$ near $x=a$ and $u(a) = u'(a) = \cdots = u^{(n-2)}(a) =0$, from the Taylor expansion of $u$ at $a$ we know that $u^{(n-1)}(a) $ cannot be negative.  If it is positive, we are done.  Otherwise, suppose $u^{(n-1)}(a)=0$, then by (\ref {eq:boundary_lemma}) we have 
\begin{eqnarray*}
u^{(n)}(a) & = & f(u(a), u'(a),..., u^{(n-1)}(a))    \\
& = & f(0, 0, ... , 0) \\
& > & 0.
\end{eqnarray*}

\vspace{.05in}
\noindent
It follows that $u(x)$ is strictly increasing near $a$ by discussions similar to those in Case 1.

\vspace{.1in}

\item Assume $u(b) = u'(b) = \cdots = u^{(n-2)}(b) = 0 $ and $u>0$ in a neighborhood of $b$.

\vspace{.05in}

\noindent
Let $s=a+b-x$, define $\hat{u}(s):=u(x)$, then
$$ \hat{u}'(s):=-u'(x), \quad \hat{u}''(s):=u''(x), \quad ... , \quad \hat{u}^{(n)}(s):=(-1)^{n}u^{(n)}(x).$$
Then $\hat{u}$ satisfies $$\hat{u}^{(n)}(s) = (-1)^{n} f \left (\hat{u}(s), -\hat{u}'(s),..., (-1)^{n-1} \hat{u}^{(n-1)}(s) \right ), $$ 
with $\hat{u}(a) = \hat{u}'(a) = \cdots = \hat{u}^{(n-2)}(a) = 0 $ and $\hat{u}>0$ in a neighborhood of $a$.

\vspace{.05in}

\noindent
Then by the result in (i) we know that either
$$\hat{u}^{(n-1)}(a)>0$$ or $$\hat{u}^{(n-1)}(a)=0, \quad \hat{u}^{(n)}(a)>0.$$

\noindent
In either case, $\hat{u}$ is strictly increasing near $a$.

\vspace{.05in}

\noindent
Therefore, either
$$(-1)^{n-1}u^{(n-1)}(b) > 0$$ or $$u^{(n-1)}(b)=0, \quad (-1)^n u^{(n)}(b) > 0.$$

\noindent
In either case, $u$ is strictly decreasing near $b$.

\end{enumerate}
 
\qed

\vspace{.1in}

In this theorem, it is necessary to assume that $f$ is Lipschitz.  For example, in Equation (\ref{eq:Holder}) the function $f(z_1, ... , z_{n})=z_1^{\alpha}$ is only H\"older continuous, but not Lipschitz continuous.  The solution 

$$u = \begin{cases}(-1)^{n-1}\lambda_n (-x)^{\frac{ n}{1-\alpha}},  & x <0 \\
     \qquad -  ~\lambda_n x^{\frac{n}{1-\alpha}},  & x \geq 0
      \end{cases} $$
satisfies $u^{(n-1)}(0)=u^{(n)}(0)=0$, so the theorem does not hold in this case.

\vspace{.2in}
\appendix \section {An Alternative Proof of Theorem \ref{thm:Hopf_k=2}}
\label{sec:appendix}

Here we give an alternative and elegant proof of Theorem \ref{thm:Hopf_k=2} suggested by the referee of an earlier manuscript.  This proof was inspired by \cite{BNV}.

\pf Let $$h_i=\sin \left (  \frac{\pi}{2} + \frac{\pi}{9} \cdot \frac{x-y_i}{x_i-a}\right )$$ where $y_i=\frac{x_i+a}{2}$.  Then
$$0 < \sin \left (  \frac{\pi}{2} - \frac{\pi}{9} \right ) \leq h_i \leq  1 < \infty \hspace{.2in} \text{on} \,\, [a, x_i]$$

\begin{eqnarray*}
h'_i(x) & = & \frac{\pi}{9} \cdot \frac{1}{x_i-a} \cos \left (  \frac{\pi}{2} + \frac{\pi}{9} \cdot \frac{x-y_i}{x_i-a}\right ) \\
h''_i(x) & = & - \left (  \frac{\pi}{9} \cdot \frac{1}{x_i-a} \right )^2 \sin \left (  \frac{\pi}{2} + \frac{\pi}{9} \cdot \frac{x-y_i}{x_i-a}\right )\\
& = & - \left (  \frac{\pi}{9} \cdot \frac{1}{x_i-a} \right )^2 h_i.
\end{eqnarray*}

\noindent
It follows that on $[a, x_i]$,
\begin{eqnarray*}
L[h_i] & := & h_i''(x)+a_1(x)h_i'(x)+a_0(x)h_i(x) \\
& = & - \left (  \frac{\pi}{9} \cdot \frac{1}{x_i-a} \right )^2 h_i + a_1(x) \cdot \frac{\pi}{9} \cdot \frac{1}{x_i-a} \cdot \cos \left (  \frac{\pi}{2} + \frac{\pi}{9} \cdot \frac{x-y_i}{x_i-a}\right ) + a_0(x)h_i(x) \\
& \leq & - \left (  \frac{\pi}{9} \cdot \frac{1}{x_i-a} \right )^2 \sin \left (  \frac{\pi}{2} - \frac{\pi}{9} \right ) + C \cdot \left (  \frac{\pi}{9} \cdot \frac{1}{x_i-a} \right ) + C  \\
& \to & -\infty \hspace{.3in} \text{as} \quad i \to \infty,   
\end{eqnarray*}
where the last inequality is true because $x_i \to a$ as $ i \to \infty.$

\vspace{.05in}
\noindent
Therefore when $i$ is large, $L[h_i] < 0$ on $[a, x_i]$.

\vspace{.05in}
\noindent
Define $$w_i=\frac{u}{h_i}.$$  Then
\begin{eqnarray*}
L[u] & = & L[h_i w_i]\\
& = & h_i w''_i+ \left ( 2h_i+a_1 h_i \right )w'_i + L[h_i]w_i .
\end{eqnarray*}
On $[a, x_i]$, since $L[u] \leq 0$ and $h_i >0$, we have 
$$ \tilde{L}_i[w_i]:= w''_i + \frac{2h_i+a_1 h_i}{h_i}w'_i + \frac{ L[h_i]}{h_i} w_i \,\, \leq \,\, 0. $$

\noindent
Since $L[h_i] < 0$ and $h_i >0$, the linear term coefficient $\frac{ L[h_i]}{h_i} <0$.  Thus the classical maximum principle and Hopf's lemma both apply to $\tilde{L}_i[w_i]$.  

\vspace{.05in}
\noindent
Because $w_i(a)=0$ and $w_i(x_i)>0$, by the maximum principle we have $$w_i(x) >0 \hspace{.2in} \text{in} \,\, (a, x_i).$$
Then by the Hopf's lemma $$w_i'(a)>0.$$  

\noindent
Finally, since $  w_i(a)=\frac{u(a)}{h_i(a)} =0 $ and $h_i(a)>0$, we obtain $$u'(a) \,\, = \,\, w_i'(a)h_i(a)+w_i(a)h_i'(a) \,\, > 0.$$

\noindent
This proves Theorem  \ref{thm:Hopf_k=2}.

\qed

\vspace{.3in}

\noindent
\textbf{Acknowledgment}: 

\vspace{.05in}
\noindent
The research of Yifei Pan was partially supported by the Senior Summer Research Grant at Indiana University - Purdue University Fort Wayne.  

\vspace{.05in}
\noindent
The research of Yu Yan was partially supported by the Ferne and Audry Hammel Research Grant at Huntington University.

\vspace{.2in}

\bibliographystyle{plain}

\bibliography{thesis}

\end{document}